\documentclass[12pt]{article}
\usepackage{amsmath,amsfonts,amssymb,amsthm}

\usepackage[hmargin= 2cm,vmargin=1cm]{geometry}

\def\a{\alpha}

\def\g{\gamma}

\def\o{\omega}

\def\CC{{\mathbb C}}
\def\FF{{\mathbb F}}
\def\MM{{\mathbb M}}

\def\SS{{\mathbb S}}

\def\ZZ{{\mathbb Z}}

\def\la{\langle}
\def\ra{\rangle}

\def\vac{{\bf 1}}
\def\vnat{{V^{\natural}}}

\def\ip#1#2{\langle #1 , #2 \rangle}

\begin{document}
\newtheorem{thm}{Theorem}[section]
\theoremstyle{definition}
\newtheorem{de}[thm]{Definition}

\newtheorem{nota}[thm]{Notation}
\newtheorem{ex}[thm]{Example}
\newtheorem{proc}[thm]{Procedure}

\theoremstyle{remark}
\newtheorem{rem}[thm]{Remark}

\begin{center}
{\Large  \bf  
Integral forms 
in vertex operator algebras, a survey}

{\it  Lecture at Ischia Group Theory Meeting, 
20 March, 2018}

\bigskip

Robert L. Griess Jr.

Department of Mathematics,

University of
Michigan,

Ann Arbor, MI 48109-1043  USA

{\tt rlg@umich.edu}

\bigskip 
\end{center} 


\begin{abstract}

We give a brief survey of recent work on integral forms in vertex operator algebras (VOAs).

\end{abstract} 

\section{Introduction}


The definition of VOA is too long to give here.   We refer the reader to a standard reference for VOA theory: \cite{FLM}, definition p.244. 
Short version: $V=\bigoplus_{n\ge 0} V_n$, a graded vector space in characteristic 0 with each $dim(V_i)$ finite;  vacuum element $\vac \in V_0$,  Virasoro element $\o \in V_2$; a  linear monomorphism $Y: V\rightarrow End(V)[[z^{-1}, z]]$, written $Y(v,z)=\sum_{n\in \ZZ} v_n z^{-n-1}$.   

For each $n\in \ZZ$, there is a product $a, b \mapsto a_n b$ (meaning the endomorphism $a_n$ applied to $b$), giving a ring $(V, n^{th})$.   

A vertex algebra (VA) is a generalization of VOA to graded modules over a commutative rings of scalars.   It has a vacuum element but not necessarily a Virasoro element.   Our examples will be over finite fields or the integers.

\begin{de} An integral form in an algebra in characteristic 0 is  the $\ZZ$-span of a basis which is closed under the product.  
\end{de} 

\begin{ex} (1) $Mat_{n\times n}(\ZZ) \le Mat_{n\times n}(\CC )$;  (2) In a simple finite dimensional complex Lie algebra, the $\ZZ$-span of a Chevalley basis is an integral form.
\end{ex}

\begin{ex} For an even integral lattice $L$, there is a lattice type VOA $V_L$ which, as a graded vector space,  has shape  
$\SS (\hat { H}) \otimes \CC [L]$.   Here,   $\SS $ means symmetric algebra of a vector space,   
$H:=\CC \otimes_{\ZZ} L$ and 
${\hat  H}$ 
means $H_1 \oplus H_2 \oplus \dots $ where $H_k$ is a copy of $H$ declared to have degree $k$.   
Finally $\CC [L]$ is the group algebra of the abelian group $L$, with basis $e^{\a}$, for $\a \in L$.
\end{ex}

\begin{ex} The Moonshine VOA $\vnat$ is a twisted version of $V_{\Lambda}$, where $\Lambda$ is the Leech lattice (rank 24, determinant 1, minimum norm 4).  It has $Aut(\vnat )\cong \MM$, the Monster. 
\end{ex}

\begin{de} An {\it integral form} $R$ in a vertex operator algebra $V=\bigoplus_{n\in \ZZ} V_n$ with a nondegenerate symmetric bilinear form is the $\ZZ$-span of a basis which is closed under all the VOA products and for all $n$, $R\cap V_n$ is an integral form of the vector space $V_n$; also $R$ must contain the vacuum element and a positive integer multiple of the  Virasoro element.    
\end{de}

\bigskip

So, an integral form in a VOA is a vertex algebra over the ring of integers.

%

Chongying Dong and RLG \cite{DG2012,DG2017} studied the following question.   Given a finite group $G$ in $Aut(V)$, is there an integral form in $V$ which is stable under $G$?   We have some general sufficient conditions.

Our main applications:  (1) for lattice type VOAs, ($L$ of rank $r$) there is a $G$-invariant integral form where $G$ has the form $2^{r}.O(L)$.   (One description: it is generated as a vertex algebra by the  $\ZZ$-span of the components of $Y(e^{\a},z)\vac $ for $\a \in L$.)

(2) For the Moonshine VOA $\vnat$, we proved that there is a Monster-invariant integral form.   (This is created by a kind of averaging, and is not (yet) described explicitly.)  

The recent preprint of Carnahan \cite{C} proves existence of an integral form in which every graded component has determinant 1.   This form is not given explicitly.    

\begin{rem}
We learned after our proof was written,  that in the 80s, Borcherds had asserted the existence of an integral form for lattice type VOAs.  Borcherds also observed that there is a Monster-invariant $\ZZ [\frac 12]$-form in $\vnat$ but claimed nothing about a form over $\ZZ$.   
\end{rem}   


If $J$ is an integral form in a VOA, it inherits a symmetric bilinear form from the VOA.  We say $J$ is {\it lattice integral} if $\ip xy \in \ZZ$ for all $x, y \in J$.    It is unclear when the restriction of the form to  $J$  is integral-valued (or can be multiplied by a scalar to become integral valued).  

Dong and RLG gave two sufficient conditions to prove lattice integrality.  (1) We showed that it is integral valued whenever the integral form $J$ is generated by quasi-primary vectors (in VOA theory, this means vectors annihilated by a certain operator $L(1)$).   This criterion applies to the Monster-invariant form we built earlier.    (2) We gave an averaging-type argument.

\section{Classical lattice type VOAs}   

This section represents joint work with Ching Hung Lam \cite{GL2014,GL2015}.    
Consider the case of lattice type VOA where $L$ is a root lattice of type ADE, and let $V:=V_L$.   Then $(V_1, 0^{th} )$ is a copy of the Lie algebra associated to $L$.  {\it If $J$ is our integral form, $J\cap V_1$ is the $\ZZ$-lattice spanned by a Chevalley basis!   So,  this $J$ is spanned as an abelian group  by a set of elements which generalizes ``Chevalley basis''.  } 
It turns out that a Chevalley group can be defined on $V_L$ with the standard generators $x_r(\pm 1)$ fixing the integral form.   

We can also take any commutative associative ring $R$ and form  $R\otimes J$, a vertex algebra over $R$,  called {\it the classical VA of type $L$ over $R$}.   We also get an action of the Chevalley group of type $L$ over $R$ on $R\otimes J$ as VA automorphisms.  

{\it When $R$ is a field we get all Chevalley groups of types ADE (with graph outer automorphisms) as full automorphism groups of these VAs over $R$.  
We also defined VA over $R$ for types BCGF.  This gives the 
 Steinberg variations (twisted Chevalley groups)  acting on VAs and being essentially the full automorphism groups.}
 
We would like to find  a series of VAs whose automorphism groups are essentially the Ree and Suzuki groups but have not (yet) done so.  

\bigskip

\begin{rem}
Our construction gives infinite dimensional graded modules for each Chevalley group and Steinberg variation over its field of definition.    These modules may be a good opportunity for study of representation theory  ($Ext$, indecomposables, etc.  ).    
\end{rem}

\section{The degree 2 component of a VOA}

Given a VOA $V=\oplus_{i\ge 0} V_i$, the $k$-th product gives a bilinear mao $V_i \times V_j {\longrightarrow} V_{i+j-k-1}$.   So, $V_n$ under the $n^{th}$ product is a finite dimensional algebra, denoted $(V_n, (n{-}1)^{th})$.   

In addition, 
(a) if $dim(V_0)=1$, $(V_1, 0^{th})$ is a Lie algebra; (b)    if $dim(V_0)=1$ and $dim(V_1)=0$, then $(V_2, 1^{st})$ is a commutative algebra with a symmetric, associative form $(ab,c)=(a,bc)$.    Algebras as in (b) are sometimes called {\it Griess algebras.   }

There are many examples of algebras (b) with finite automorphism groups, e.g., \cite{S2004,S2006}, \cite{DGR, DG2005}.  The 196884-dimension algebra used to construct the Monster occurs this way in the Moonshine VOA $V^{\sharp}$.   


\bigskip

Now suppose $dim(V_0)=1$ and that $e\in V_2$ is a conformal vector of central charge $\frac 12$ and that the subVOA generated by $e$ is simple.  Miyamoto showed that $e$ gives $t_e \in Aut(V)$ of order 1 or 2 (called a {\it Miyamoto involution}).   

\bigskip

In the special case of a dihedral VOA (generated by a pair of such conformal vectors $e, f$), the degree 2 algebra has integral forms.  Those which are maximal  integral forms and invariant under the dihedral group $\la t_e, t_f \ra$ were classified in the thesis of Greg Simon (U Michigan, 2016).    For the nine types of dihedral VOAs (classified by Sakuma \cite{Sak}; they correspond to nodes of the {\it extended $E_8$-diagram}, displayed below), there is just one maximal invariant form in all cases but $2A$, in which case there are three.

{\Large 
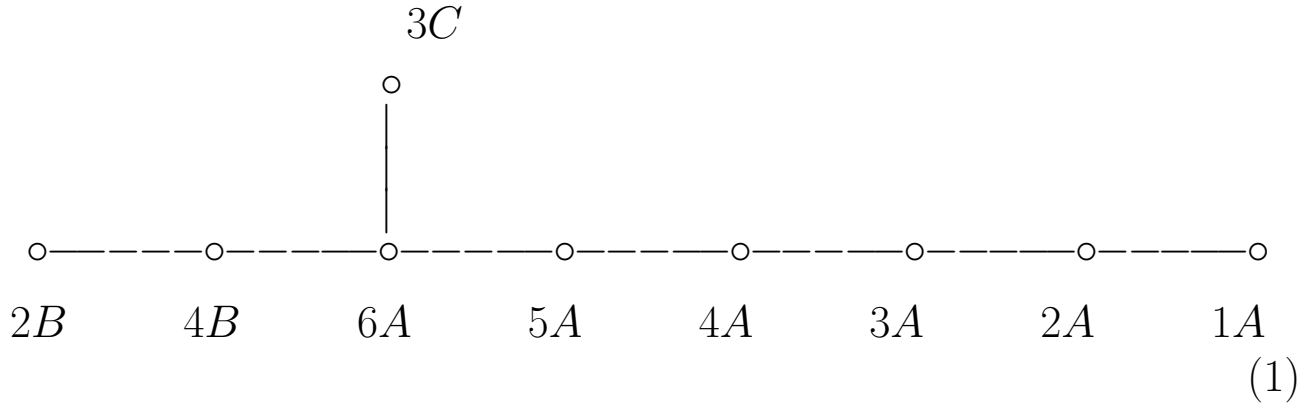
\begin{equation}\label{eq:1.1}
\begin{array}{l}
  \hspace{150pt} 3C\\
  \hspace{140pt}\circ \vspace{-6.2pt}\\
   \hspace{140pt}| \vspace{-6.1pt}\\
 \hspace{140pt}| \vspace{-6.1pt}\\
 \hspace{140pt}| \vspace{-6.2pt}\\
  \hspace{6pt} \circ\hspace{-5pt}-\hspace{-7pt}-\hspace{-5pt}-
  \hspace{-5pt}-\hspace{-5pt}-\hspace{-5pt}\circ\hspace{-5pt}-
  \hspace{-5pt}-\hspace{-5pt}-\hspace{-6pt}-\hspace{-7pt}-\hspace{-5pt}
  \circ \hspace{-5pt}-\hspace{-5.5pt}-\hspace{-5pt}-\hspace{-5pt}-
  \hspace{-7pt}-\hspace{-5pt}\circ\hspace{-5pt}-\hspace{-5.5pt}-
  \hspace{-5pt}-\hspace{-5pt}-\hspace{-7pt}-\hspace{-5pt}\circ
  \hspace{-5pt}-\hspace{-6pt}-\hspace{-5pt}-\hspace{-5pt}-
  \hspace{-7pt}-\hspace{-5pt}\circ\hspace{-5pt}-\hspace{-5pt}-
  \hspace{-6pt}-\hspace{-6pt}-\hspace{-7pt}-\hspace{-5pt}\circ
  \hspace{-5pt}-\hspace{-5pt}-
  \hspace{-6pt}-\hspace{-6pt}-\hspace{-7pt}-\hspace{-5pt}\circ
  \vspace{-6.2pt}\\
  \vspace{-10pt} \\
  2B\hspace{44pt} 4B\hspace{44  pt} 6A\hspace{44pt} 5A\hspace{44pt}
  4A\hspace{44pt} 3A\hspace{44pt} 2A\hspace{44pt} 1A\\
\end{array}
\end{equation}
}

This classification of maximal invariant forms in $V_2$ does not (yet) extend to invariant forms in the entire dihedral VOA.

\section{Modular Moonshine of Borcherds and Ryba}

Borcherds and Ryba wrote several articles \cite{BR,R} about Modular Moonshine (positive characteristic) which imitated the story of the Monster and the graded representation $\vnat$ and modular forms, but for smaller sporadic groups.   

They discussed an interesting case.  In $\MM$, take $g$ a $3C$-element; then $C(g)=\la g \ra \times S$, where $S\cong F_3$ a sporadic simple group of order $2^{15}3^{10}5^3 7^2 13\cdot 19\cdot 31 $ (Thompson's group).    

Borcherds and Ryba used $K$, Borcherds's  $\ZZ [\frac 12]$-form in $\vnat$, then considered  its 0th Tate cohomology group $${\hat H}^0(\la g \ra, K):=K^g/(1+g+g^2)K.$$  This inherits  structure to make a VA over $\FF_3$.   It looked like the classical $E_8$ type    VA over $\FF_3$, but  nonzero terms occur only in degrees 0, 3, 6 , . . .  and have respective dimensions $1, 248, \dots $, just like for the genuine $E_8$ VA in degrees 0, 1, 2, . . ..     There was no obvious  isomorphism (which triples degree of the grading) between these two VAs over $\FF_3$.    

To prove existence of an isomorphism, Lam and RLG adapted a covering idea of Frohardt-Griess \cite{FG1992}, which is illustrated in the following example.

\begin{ex} $F$ algebraically closed field of char.3.   The Lie algebra $a_2(F)$ has a 1-dimensional central ideal, $Z$. 
While $Aut(a_2(F))$ is  $PGL(2,F){:}2$, $Aut(a_2(F)/Z)\cong G_2(F)$.  Our proof takes the Lie algebra $d:=d_4(F)$ and graph automorphism $\g$ of order 3, then considers 
$$ 0 < (1+\g +\g^2)d < d^{\g} < d, \text { \ dimensions 0, 7, 14, 28.} $$   The group $G_2(F)$ acts on each subobject and on the 7-dimensional quotient Lie algebra $d^{\g}/(1+\g +\g^2)d$.   One can see inside $d$ (look at the long roots) a copy of $a_2(F)$ which maps onto $d^{\g}/(1+\g +\g^2)d$; the image is isomorphic to $a_2(F)/Z$.   Therefore $Aut(a_2(F)/Z)$ contains a copy of $G_2(F)$.   This containment is equality.    
\end{ex} 

Lam and RLG  took $M$, the standard integral form for $V_{E_8}$, and a sublattice $M'$ of $K$ (integral form in $\vnat$) which ``covered'' the  Tate cohomology group $K^g/(1+g+g^2)K$.   ({\it Roughly,}   $M'$ is the $\ZZ$-span of the image of $M$ under a map suggested by $x\mapsto x\otimes x \otimes x$ for $x\in EE_8\cong \sqrt 2 E_8$;  think of the containment of lattices $EE_8 \perp EE_8 \perp EE_8 \le \Lambda$, the Leech lattice).  This led to an isomorphism.   

An application of the Borcherds-Ryba theory is a new proof that the group $F_3$ of Thompson embeds in $E_8(3)$ (first proof 1974, by Thompson and P. Smith, used a study of Dempwolff decompositions and computer work).   This VA viewpoint gives a nontrivial homomorphism  of $C(g)/\la g \ra \cong F_3$ into the group $E_8(3)$ without knowing much about the structure of $C(g)$.    

\bigskip 
 
My web site contains files of certain articles in the reference list.

\end{document}